 \documentclass[11pt]{article}

\usepackage{amsmath}
\usepackage{amssymb}
\usepackage{amsthm}
\usepackage{graphicx}
\usepackage{color}
\usepackage{mathrsfs}
\newcommand{\be}{\begin{equation}}
\newcommand{\ee}{\end{equation}}
\newcommand{\bes}{\begin{equation*}}
\newcommand{\ees}{\end{equation*}}
\newcommand{\beqn}{\begin{eqnarray}}
\newcommand{\eeqn}{\end{eqnarray}}
\newcommand{\beqns}{\begin{eqnarray*}}
\newcommand{\eeqns}{\end{eqnarray*}}

\newcommand{\lfi}{\left\{}
\newcommand{\rfi}{\right\}}

\newcommand{\Del}{\Delta}

\newcommand{\ga}{\gamma}

\newcommand{\sig}{\sigma}

\newcommand{\MISE}{\mbox{MISE}}

\newcommand{\Var}{\mbox{Var}}

\newtheorem{theorem}{Theorem}
\newtheorem{lemma}{Lemma}

\newcommand{\etal}{{\it et al.\ }}

\long\def\ignore#1{}

\begin{document}

\title{{\bf Density Deconvolution  with Small Berkson Errors}}
\author {Ramchandra Rimal and Marianna Pensky  \\
Department  of Mathematics\\ 
University of Central Florida\\
 Orlando, FL 32816-1364 }
 
\date{}

\maketitle

\vspace{10mm}

% \noindent
% {\bf RUNNING HEAD:}\ {\bf Deconvolution  with Small Berkson Errors }

%\vspace{7cm} 
%\noindent
% This work was supported in part by National Science Foundation (NSF),   
% grants     DMS-1407475 and DMS-1712977.   

\begin{abstract}
The present paper studies density deconvolution in the presence of small Berkson errors, 
in particular, when the variances of the errors tend to zero as the sample size grows.
It is known that when the Berkson errors are present, in some cases, the unknown density estimator can be obtain by 
simple averaging without using kernels. However, this may not be the case when Berkson errors are asymptotically small. 
By treating the former case as a kernel estimator with the zero bandwidth, we  obtain the optimal expressions
for the bandwidth. We  show that the density of  Berkson errors acts as a regularizer, so that the kernel estimator 
is unnecessary when the variance of Berkson errors lies above some threshold that depends on the on the shapes 
of the densities in the model and the number of observations.
\\ 


\end{abstract}

\newpage

%%%%%%%%%%%%%%%%%%%%%%%%%%%%%%%%%%%%%%%%%%%%%%%%%%%%%%%%%%%%%%%%%%%%%%%%%%%%%%%%%%%%%

\section{Introduction}
\label{sec:intro}
\setcounter{equation}{0}

In many real life problems one is interested in distribution of a certain variable 
which can be observed only indirectly. Mathematically, this leads to a density deconvolution problem 
where one needs to estimate the pdf of a variable $X$ on the basis of observations of a surrogate variable $Y = X + \xi$
where the pdf $f_\xi$ of $\xi$ is known. The real life applications of this model arise in econometrics, astronomy,
biometrics, medical statistics, image reconstruction (see, e.g., Bovy \etal (2011), Robinson (1999), Wason \etal (1984), 
and also Carroll \etal (2006) and Meister (2009) and references therein). 
Density deconvolution problem was extensively studied in the last thirty years (see, e.g., 
Carroll \etal (2009), Comte and  Kappus  (2015), Goldenshluger (1999),  Lacour and Comte  (2011), among others,  
and Meister (2009) and references therein).

However,   Berkson (1950) argued that in many situations it is more appropriate to treat the true unobserved variable 
as being contaminated with an error itself and search for the distribution of $W = X + \eta$ where $\eta$ is the so-called 
Berkson error with a known pdf $f_\eta$. Here, $X$, $\xi$ and $\eta$ are assumed to be independent. The objective is to  estimate  
the pdf $f_W$ of $W$ on the basis of i.i.d. observations 
\be \label{eq:deconv}
Y_i=X_i + \xi_i,\ \ i=1, \cdots, n, 
\ee
where $X_i$ and $\xi _i$ are   i.i.d with, respectively,  the pdfs $f_X$ which is unknown and   $f_\xi$ which is known.
Density $f_\xi$   is called the error (or the blurring) density.

Estimation with Berkson errors occurs in a variety of statistics fields such as 
analysis of chemical, nutritional, economics or astronomical data (see, e.g., Kim \etal (2016),
Long \etal (2016), Robinson (1999), Wang (2003), Wason \etal (1984) among others). 
For example, in occupational medicine, an important problem is the assessment of
the health hazard of specific harmful substances in a working area. A 
 modeling approach usually assumes that there is a threshold concentration,
called the threshold limiting value (TLV), under which there is no risk
due to the substance. Estimating the TLV is of particular interest in the
industrial workplace. The classical errors in this model come from the measures of
dust concentration in factories, while the Berkson errors come from the
usual occupational epidemiology construct, wherein no direct measures
of dust exposure are taken on individuals, but instead plant records of
where they worked and for how long are used to impute some version
of dust exposure (see Carroll \etal (2006)). In economics, the household income is usually not precisely collected due to
the survey design or data sensitivity. It was described by Kim \etal (2016) (see also Geng and Koul (2018)) that
when the income data were collected by asking individuals which salary range categories
they belong to, then the midpoint of the range interval was used in analysis. 
In this case, it is wise to assume that the true income fluctuates around 
the midpoint observation subject to errors.

Estimation with Berkson errors was studied by Carroll \etal (2009), Delaigle (2007, 2008),
Du \etal (2001), Geng and Koul  (2018),  Wang (2003, 2004) among others.  
It is well known that the presence of Berkson errors improves precision of estimation 
of the density function $f_W$. For example, Delaigle (2007, 2008) who studied estimation with Berkson errors 
noted that in the cases when the pdf $f_\eta$ of Berkson errors has higher degree of smoothness than 
the error density  $f_\xi$,   one can   obtain  estimators of $f_W$ with the parametric convergence rate.

However, in majority of practical situations, the Berkson errors are small.
Hence, the question arises whether small Berkson errors improve the estimation accuracy and how much.
A similar inquiry has been recently carried out by Long  \etal  (2016)  who considered a somewhat  
different setting. In particular, they studied a $p$-dimensional version of the problem where 
variable $X$ is directly observed and the objective is  estimation of the pdf $f_W$ of $W = X + \eta$ on the basis of observations 
$X_1, \cdots, X_n$ where the pdf $f_\eta$ of $\eta$ is known and variable $\eta$ is small. In this formulation, the  pdf $f_W$ can be written as 
$$
f_W(x) = \int_{\mathbb{R}^p} f_X(x-z)f_{\eta}(z) dz
$$
and can be estimated  by 
\be \label{eq:Long1}
\hat{f}_W(x) =  n^{-1} \sum_{i=1}^{n}f_{\eta}(x - X_i)
\ee
with the parametric error rate of $C n^{-1}$. However,   if $\Var(\eta) =\sigma^2$  is small, 
this rate becomes $C(\sigma) n^{-1}$ 
 where $C(\sigma) \to \infty$ when $\sigma \to 0$, so the error  of the estimator \eqref{eq:Long1}  may be very high.

To resolve this difficulty, in addition to estimator \eqref{eq:Long1}, Long {\it et al.} (2016) proposed two 
alternative kernel estimators where the bandwidths of the kernels are chosen as $h = h_W$ or $h = h_X$, 
so to minimize the error of the estimator of $f_W$ in the first case and the error of the estimator of $f_X$  in the second case.  
Subsequently, the authors studied all three estimators by simulations and concluded 
that overall the kernel estimator with $h = h_W$ outperforms the  remaining two. 
When the error variance $\sig$ is small, the estimator \eqref{eq:Long1} leads to sub-optimal error rates.
On the other hand, the choice of $h = h_X$ leads to oversmoothing, especially, when 
the error variance is large. The authors do not provide a comprehensive theoretical study of the bandwidth selection in a general case.
In particular, their  rule-of-thumb recipe is based on the case where $f_X$ is a Gaussian density. 
In particular,  Long  \etal  (2016)  did not investigate when estimator \eqref{eq:Long1}  that corresponds to the bandwidth $h=0$  is preferable 
and suggested that it is always suboptimal.

The objective of the present paper is to study the situation 
where both the blurring and the Berkson errors are present
and, in addition, the  Berkson errors  $\eta_i$, $i=1, \cdots, n$, are   small.  
To quantify this phenomenon, we assume that   the pdf  $f_\eta$ is of the form 
\be \label{eq:f_eta}
 f_\eta (x)= \sigma^{-1} g\left(\sigma^{-1}  x\right), 
\ee
where $\sigma$ is small, specifically, $\sigma = \sigma_n \to 0$ as $n \to \infty$
while the variable $X$ has a non-asymptotic scale. 
Specifically, we shall provide a full theoretical study of the bandwidth selection 
in a density deconvolution with Berkson errors.

The setting of Long {\it et al.} (2016) corresponds to the multivariate version of the problem in this paper 
where $\xi_i = 0$ and  $f^*_\xi =1$. 
We provide   full theoretical treatment of the problem. In particular, we prove that 
one should always choose the bandwidth to minimize the error  of the estimator of $f_W$, 
but in some cases this optimal bandwidth can be   zero if $\sigma$ lies above some threshold 
that depends on the shapes of the densities  $f_\xi$, $f_X$ and $g$ and the sample size. In the particular case studied by Long \etal (2016), 
the latter situation would lead to the estimator of the form \eqref{eq:Long1}.

Since the setting \eqref{eq:f_eta} leads to three asymptotic parameters, $n$, $\sigma$ and $h$,
in order to keep the paper clear and readable,  we consider a one-dimensional version of the problem. 
Extensions of our results to the situation of multivariate densities is a matter of future work.

In what follows, we are using the following notations. For any function $f$, $f^*$ denotes its Fourier transform.
If $f$ is a pdf, then $f^*$ is the characteristic function of $f$. 
We use the symbol $C$ for a generic positive constant, 
which   takes different values at different places and is independent of $n$.
Also, for any positive functions $a(n)$ and $b(n)$, we write $a(n) \asymp b(n)$ if the ratio 
$a(n)/b(n)$ is bounded above and below by finite positive  constants independent of $n$, and 
  $a(n) \lesssim b(n)$ if the ratio $a(n)/b(n)$ is bounded above  by finite positive  constants independent of $n$.

The rest of the paper is organized as follows. Section~\ref{sec:construction} presents an  
estimator of $f_W$  in the case of the small Berkson errors. Section~\ref{sec:est_error} provides an expression for the
error of this estimator and also derives the optimal value of the bandwidth that depends on the shapes of 
the densities  in the model  and on the values of parameters $n$ and $\sigma$. 
Section~\ref{sec:discuss} is devoted to the discussion of the results of the paper.
The proofs of all statements can be found in Section~\ref{sec:proofs}.

%%%%%%%%%%%%%%%%%%%%%%%%%%%%%%%%%%%%%%%%%%%%%%%%%%%%%%%%%%%%%%%%%%%%%%%%%%%%%%%%%%%%%%%%%%%%%%%%%%%%%%%%%%%%

\section{Construction of the deconvolution estimator} 
\label{sec:construction}
\setcounter{equation}{0}

Since \eqref{eq:deconv} and $W = X + \eta$  imply  that
\be \label {Fequation1}
f_Y^* (w) = f_X ^*  (w)  f_\xi^*(w), \quad
f_W^* (w) = f_X ^*  (w)  f_\eta^*(w) 
\ee
and also, due to \eqref{eq:f_eta},  $f_\eta^* (w)= g^* (\sigma w)$, one obtains
$$
f_W^* (w) = f_X ^*  (w) g^* (\sigma w) =\frac{ f_Y ^*  (w) g^* (\sigma w)} {f_\xi ^*  (w) }
$$
Note that the unbiased estimator of $f_Y ^*  (w) $ is given by the empirical characteristic function 
\be \label{charact}
 \hat{f}_Y ^*  (w)= n^{-1} \sum_{j=1}^{n} \exp(iwY_j).
\ee
If $g^* (\sigma w)/f_\xi ^*  (w)$ is square integrable, i.e.
\be \label{rho_delta}
\rho^2(\sigma) =\int_{-\infty}^{\infty} \left| \frac{  g^* (\sigma w)} {f_\xi ^*  (w) }\right| ^2 dw < \infty, 
\ee
then the inverse Fourier transform of $f_Y ^*  (w) g^* (\sigma w)/f_\xi ^*  (w)$
exists and $f_W(x) $ can be estimated by 
\be \label{estimator}
\hat{f}_W (x) =\frac{1}{2\pi} \int _{-\infty}^{\infty} \exp(-iwx)\ \frac{ \hat{f}_Y ^*  (w) g^* (\sigma w)} {f_\xi ^*  (w) } dw
\ee
If  $g^* (\sigma w)/f_\xi ^*  (w)$ is not square integrable, one needs to obtain a kernel estimator of $f_W$.  
Construct   approximations $f_{W,h}$ and $f_{W,h}^*$ of $f_{W}$  and  $f_{W}^*$, respectively,   
\be\label{Fequation 4}
f_{W,h} (x) = \int_{-\infty}^{\infty} \frac{1}{h} K\left(\frac{x-w}{h}\right) f_W(w) dw,\quad
f_{W,h}^* (s) =   K^*(sh) \frac{ f_Y ^*  (s) g^* (\sigma s)} {f_\xi ^*  (s) }
\ee
and arrive at the estimator $\hat{f}_{W,h}^* (s)$ of $f_{W,h}^* (s)$  of the form 
\bes
\hat{f}_{W,h}^* (s) = K^*(sh)  \hat{f}_Y ^*  (s) g^* (\sigma s)/f_\xi ^*  (s)
\ees
where $\hat{f}_Y ^*$ is defined in  \eqref{charact}.

Consider the kernel function $K(x) =  \sin(x)/(\pi x)$,  
so that  $K^* (s) =I(|s|   \leq 1),$  where $I(A)$ denotes the indicator function of a set $A$. 
Since $K^* (s)$ is bounded and compactly supported, the inverse Fourier transform of  $\hat{f}_{W,h}^*$ always exists 
and 
\be  \label{ker_est}
\hat{f}_{W,h}(x) = \frac{1}{2\pi}\, \int_{-\infty}^{\infty}  \exp(-ixs) \    
\frac{ \hat{f}_Y ^*  (s) K^* (sh)  g^* (\sigma s)}{f_\xi ^*  (s)}  \,  ds
% = \frac{1}{2\pi}\int_{-\infty}^{\infty}  e^{-i s w}\,   
%\frac{ \hat{f}_Y ^*  (s)  g^* (\sigma s)}{f_\xi ^*  (s)}\, I(|s| < h^{-1}) \,  ds.
\ee
Note that   $\hat{f}_{W,0} (x) =\hat{f}_W(x)$.

In order to obtain an expression for the bandwidth $h$  we  introduce   the following assumptions:
\\

\noindent
{\bf (A1)}\ 
There exists positive numbers $c_\xi$ and $C_\xi$  and nonnegative numbers $a, b$ and $d$ such that for any~$s$
\begin{equation} \label{Ass_A1}
 c_\xi (s^2 +1)^{-\frac{a}{2}  }    \exp(-d | s| ^b) \leq  |f^*_{\xi} (s)|  \leq C_\xi (s^2 +1)^{-\frac{a}{2}  }\exp(-d | s| ^b)  
\end{equation} 
where $b=0$ iff $d=0$ and $a>0$ whenever $d=0$.   
\\
 
\noindent
{\bf (A2)}\ 
There exists positive numbers $c_g$ and  $C_g$ and nonnegative numbers $\alpha,\beta$ and $\gamma$  such that for any~$s$ 
\begin{equation} \label{Ass_A2}
c_g (s^2+1)^{-\frac{\alpha}{2}  }    \exp(-\gamma |s| ^\beta) \leq   |g^* (s)|  \leq C_g (s^2+1)^{-\frac{\alpha}{2}}   \exp(-\gamma | s| ^\beta)  ,\\  
\end{equation}
where $\beta=0$ iff $\gamma=0$ and $\alpha>0$ whenever $\gamma=0$.\\

\noindent
{\bf (A3)}\ 
$f_{X} (s)$   belongs to the Sobolev ball 
\begin{equation} \label{Ass_A3}
\mathcal{S}(k,B)=\{f: \int_{-\infty} ^\infty\ |f^*_X(s)|^2( s^2+1)^k\, ds \leq B^2\}.  
\end{equation} 
\\
  
Also, since density deconvolution with  Berkson errors of relatively large size  has been fairly well studied,
below we   only  study the case where  $\sigma$ is small, in particular, if $\ga>0$, $d>0$, one has
\be \label{sig-ineq}
\sigma < 0.5\,  (d/\gamma)^{1/b}.
\ee

%%%%%%%%%%%%%%%%%%%%%%%%%%%%%%%%%%%%%%%%%%%%%%%%%%%%%%%%%%%%%%%%%%%%%%%%%%%%%%%%%%%%%%%%%%%%%%%%%%%%%%%%%%%%%%%%%%%%%%%%%%%%%%%%%%%%%%%

\section { Estimation error}
\label{sec:est_error}
\setcounter{equation}{0}

We characterize the precision of the  estimator $\hat{f}_{W,h}$ of $f_W$   by its Mean Integrated Squared Error (MISE)
\bes 
\MISE  (\hat{f}_{W,h}, f_W) = E \int_{-\infty}^{\infty} | \hat{f}_{W,h}(x) - f_W(x)| ^2 dx.
\ees
Since, under Assumptions \eqref{Ass_A1}--\eqref{Ass_A3}, both $ \hat{f}^*_{W,h}$ and $f^*_W$ are square integrable, 
by  the Plancherel theorem, derive that
\begin{align*}
 \MISE  (\hat{f}_{W,h}, f_W)  
=  \frac{1}{2\pi}\,  E   \int_{-\infty}^{\infty}   \frac{  |g^* (\sigma s)|^2} { |f_\xi ^*  (s)|^2}\,   | K^*(sh) \hat{f}_Y ^*  (s)  -  f_Y ^*  (s)|^2\, ds
\end{align*} 
Therefore,  
\be \label{mise}
\MISE  (\hat{f}_{W,h}, f_W)  = R_1 (\hat{f}_{W,h}, f_W) +  n^{-1}\,  R_2 (\hat{f}_{W,h}, f_W)
\ee 
where  
\bes
R_1 (\hat{f}_{W,h}, f_W)
  =   \| E \hat{f}_{W,h} -f_W\|^2 
 =   \frac{1}{2\pi}\, \int_{-\infty}^{\infty} \,   | g^* (\sigma s)|^2 \,  |f_X ^*  (s) |^2 \, I(|s| > h^{-1})\, ds 
\ees
  is the integrated squared  bias of the estimator $\hat{f}_{W,h}$ and  
  \bes
   R_2 (\hat{f}_{W,h}, f_W) 
 =  E \| \hat{f}_{W,h}- E \hat{f}_{W,h} \|^2 
 %    &\leq \frac{1}{2n\pi} \int_{-\infty}^{\infty}   \frac{  | g^* (\sigma s) | ^2} { | f_\xi ^*  (s)| ^2}    |  K^*(sh)    | ^2ds \\
 \leq \frac{1}{2\pi} \int_{-\infty}^{\infty}   \frac{|g^* (\sigma s)|^2} {|f_\xi ^*  (s)|^2}\  I(|s| < h^{-1}) \, ds
   \ees
   is the integrated variance. 
%%%%%%%%%%%%%%%%%%%%%%%%%%%%%%%%%%%%%%%%%%%%%%%%%%%%%%%%%%%%%%%%%%%%%%%%%%%%%%%%%%%%%%%%%%%%%%%%%%%%%%%%%%%%%%%%%%%%%%%%%%%%%%%%%%

\begin{table} [h]
\begin{center}
\begin{tabular}{|l|c| }
\hline
Case &  $\Delta_2$ \\
\hline
I)\ $b = \beta =0, \alpha > a + \frac{1}{2} $,  & $  \min \left( h^{-(2a +1)},  \sigma^{-(2a+1)}\right)$ \\
 & \\
\hline
%%%%%%%%%%%%%%%%%%%%%%%%%%%%%%%%
II) $b = \beta =0 , \alpha = a + \frac{1}{2}$  & $  \min \left( h^{-(2a +1)},  \sigma^{-(2a+1)}\right)\max \left\{\ln\left({\frac{\sigma}{h} }\right),1\right\}$ \\
&\\
\hline
%%%%%%%%%%%%%%%%%%%%%%%%%%%%%%%%
III)  
$b = \beta =0, \alpha < a + \frac{1}{2}$,  & $h^{-(2a+1)}\,  \min\left\{\left(   \frac{h}{\sigma}\right)^{2\alpha}  , 1  \right\}   $ \\
&\\
\hline
%%%%%%%%%%%%%%%%%%%%%%%%%%%%%%%%%%%%%%%%%
IV) 
$b =0, \beta > 0$ & $  \min \left( h^{-(2a +1)},  \sigma^{-(2a+1)}\right)$ \\
& \\ 
\hline
%%%%%%%%%%%%%%%%%%%%%
V)   
$\beta > b > 0 , h > \left(\frac {\gamma \beta}{db}{\sigma^{\beta}}\right)^{\frac{1}{\beta-b}} $ 
& $ h^{-(2a+1)+ b}{\exp (2dh^{-b}}) \min\left\{\left(   \frac{h}{\sigma}\right)^{2\alpha}, 1  \right\} $ \\
\ \ \ \  $\beta > b > 0 , h < \left(\frac {\gamma \beta}{db}{\sigma^{\beta}}\right)^{\frac{1}{\beta-b}} $ 
& $  \exp\left(\kappa   \sigma^{ -\frac{\beta b}{\beta - b}}\right)\sigma^{\frac{\beta}{\beta -b}.\frac{b-2}{2} -2\alpha}$ \\
%%%%%%%%%%%%%%%%%%%%%%%
\hline
VI)
$b = \beta > 0$ &   $h^{-(2a+1)+ b}\exp(2h^{-b}(d- \gamma \sigma^{b})) \min\left\{\left(   \frac{h}{\sigma}\right)^{2\alpha}  , 1  \right\}$ \\
 & \\
%%%%%%%%%%%%%%%%%%%%%%%%%%
\hline 
  VII)
  $b > 0,\beta = 0$ & $ h^{-(2a+1)+ b}{\exp (2dh^{-b}}) \min\left\{\left(   \frac{h}{\sigma}\right)^{2\alpha}  , 1  \right\}$  \\
 &  \\
 %%%%%%%%%%%%%%%%%%%%% 
 \hline
 VIII)
  $b > \beta >  0$ &  $ h^{-(2a+1)+ b}{\exp (2dh^{-b}}) \min\left\{\left(   \frac{h}{\sigma}\right)^{2\alpha}  , 1  \right\}$  \\
 &  \\
 \hline
\end{tabular}
%%%%%%%%%%%%%%%%%%%%%%%%%%%%%%%%%%%%%%%%%%%%%%%%%%%%%%%%%%%%%%%%%%%%%%%%%%%%%%%%%%%%%%%%%%%%%%%%%%%%%%%%%%%%%%%%%%%%%%%%%%%%%%%%%%
%
\caption{The    asymptotic expressions for $\Del_2 \equiv \Del_2(\sig, h)$}
\end{center}
 \end{table}

%%%%%%%%%%%%%%%%%%%%%%%%%%%%%%%%%%%%%%%%%%%%%%%%%%%%%%%%%%%%%%%%%%%%%%%%%%%%%%%%%%%%%%%%%

We shall be interested in the maximum value of $\MISE  (\hat{f}_{W,h}, f_W) $ over all 
$f_X \in \mathcal{S}(k,B)$ where $\mathcal{S}(k,B)$ is defined in \eqref{Ass_A3}. 
In particular, we define
\be \label{eq:Del}
\Del \equiv \Del(n,\sig, h) = \max_{f_X \in \mathcal{S}(k,B)}\, \MISE  (\hat{f}_{W,h}, f_W)\quad \mbox{subject to} \quad 
f_W^* (w) = f_X ^*  (w)  f_\eta^*(w). 
\ee
It is easy to see that 
\be  \label{eq:Delta}
  \Del \leq  \Del_1 + n^{-1}\, \Del_2
\ee
\be
\label{eq:Del12}
  \Del_1   = \max_{f_X \in \mathcal{S}(k,B)}\, R_1 (\hat{f}_{W,h}, f_W),
\quad 
\Del_2   = \max_{f_X \in \mathcal{S}(k,B)}\, R_2 (\hat{f}_{W,h}, f_W)
\ee

\noindent
Then, the following statements hold.

 \begin{lemma}\label{lbias} 
Under the assumptions \eqref{Ass_A1}--\eqref{sig-ineq},  for $\Del_1$   in \eqref{eq:Del12}, one has
\be \label{bias}
\Delta_1 = \Del_1 (\sig, h)  \lesssim 
\left\{
 \begin{array}{ll}
 \sigma^{ -2\alpha} h^{2\alpha +2k} \exp\left(-2\gamma \left(\sigma/h \right)^{\beta}\right) &   \text{if}\   h < \sigma\\
 h^{2k}\	 	 	&  \text{if}\   h\geq \sigma
\end{array} \right. 
\ee
\end{lemma}

\begin{lemma}\label{variance}
If $\beta >b>0$, denote 
 \be\label{kappa}
 \kappa = \left(\frac {db}{\gamma \beta}\right)^{\frac{b}{\beta-b}} \left[ \frac{d(\beta -b)}{b}\right] > 0, \quad \beta >b >0.
 \ee 
Then, under the assumptions \eqref{Ass_A1}--\eqref{sig-ineq},  the expressions for  $\Del_2$  defined in \eqref{eq:Del12}, 
are given in  Table~1. 
 \end{lemma}

%%%%%%%%%%%%%%%%%%%%%%%%%%%%%%%%%%%%%%%%%%%%%%%%%%%%%%%%%%%%%%%%%%%%%%%%%%%%%%%%%%%%%%%%%
\noindent 
Observe that in every case, the expression for the variance depends not only on the values of 
$h$, $\sig$ and $n$ but also on their mutual relationship. 
 Also, the bias term $\Delta_1(\sigma , h) $ is an increasing function of $h$ while the variance term
$\Delta_2(\sigma , h) $  is a decreasing function of $h$, so the optimal value   $h= h_\text{opt}$ 
is such that $\Delta_1(\sigma , h) \asymp n^{-1}\, \Delta_2(\sigma, h)$.   Theorem~\ref{th:mise} below 
presents the optimal expressions $h_{opt}$ for the  bandwidth $h$ as well as the 
corresponding values for the risk $\Del(n,\sig, h_{opt})$  where $\Del(n,\sig, h)$ is defined in \eqref{eq:Del}.

%%%%%%%%%%%%%%%%%%%%%%%%%%%%%%%%%%%%%%%%%%%%%%%%%%%%%%%%%%%%%%%%%%%%%%%%%%%%%%%%%%%%%%%%%%%%%%%%%%%

\begin{theorem} \label{th:mise}
Let conditions \eqref{Ass_A1}--\eqref{sig-ineq} hold. Then, the asymptotic values of 
$$
h_{opt} = \arg\min_h \Del(n,\sig, h)
$$
and also of $\Del(n,\sig, h_{opt})$ are provided in  Table~2.
Here,
\begin{eqnarray} 
 \mu_1 = \mu_1(n) & = &  \left[\frac{1}{2d}\left( \ln n + \left(\frac{b -2a-2k -1}{b}\right)\ln\ln n\right)\right]^{-\frac{1}{b}}, \nonumber\\
& & \label{eq:mu}\\
 \mu_2 = \mu_2 (n) & = & \left[\frac{1}{2(d- \gamma \sigma^{b})}\left(\ln n + \left(\frac{b-2a -2k-1}{b}\right) \ln\ln n\right)\right]^{-\frac{1}{b}}.  
\nonumber
\end{eqnarray}
\end{theorem}

%%%%%%%%%%%%%%%%%%%%%%%%%%%%%%%%%%%%%%%%%%%%%%%%%%%%%%%%%%%%%%%%%%%%%%%%%%%%%%%%%%%%%%%%%%%%%%%%%%%

 \begin{table} [h]
\begin{center}
\begin{tabular}{|l|c|c|c| }
%\multicolumn{4}{l}{ The   optimal values of the bandwidth $h$ and the corresponding expressions for  the MISE}\\
\hline
Case &  $\Del(n, \sig, h_{opt})$  &  condition & $h_{opt}$\\
\hline
I)\ $b = \beta =0$,  & $ n^{-1}\sigma^{-(2a+1)}$ & $\sigma \geq n^{-\frac{1}{2k+2a+1}}$ & $ 0$ \\
 $\alpha > a + \frac{1}{2}$ % &  & & \\
& $n^{-\frac{2k}{2k+2a+1}}$&$ \sigma  <  n^{-\frac{1}{2k+2a+1}}$&$  n^{-\frac{1}{2k+2a+1}}$\\
\hline
%%%%%%%%%%%%%%%%%%%%%%%%%%%%%%%%
II) $b = \beta =0$  & $ n^{-1} \sigma^{-(2a+1)}\ln n$&  $\sigma  \geq  n^{-\frac{1}{2k+2a+1}}$& $ n^{-\frac{1}{2k+2\alpha}} \sigma^{\frac{2\alpha -2a-1}{ 2\alpha +2k}} $\\
 $ \alpha = a + \frac{1}{2}$ % &  & & \\ 
 & $n^{-\frac{2k}{2k+2a+1}}$ & $\sigma < n^{-\frac{1}{2k+2a+1}}$& $ n^{-\frac{1}{2k+2\alpha+1}}$\\
\hline
%%%%%%%%%%%%%%%%%%%%%%%%%%%%%%%%
III)  
$b = \beta =0$,  & $\sigma ^{-2\alpha } n ^{-\frac{2\alpha +2k}{2k +2a +1}}$ & $ \sigma  >  n^{-\frac{1}{2k+2a+1}}$ & $n^{-\frac{1}{2k+2a+1}}$\\
$ \alpha < a + \frac{1}{2}$ & $n^{-\frac{2k}{2k+2a+1}} $ & $ \sigma  \leq  n^{-\frac{1}{2k+2a+1}}$ & $n^{-\frac{1}{2k+2a+1}}$ \\
\hline
%%%%%%%%%%%%%%%%%%%%%%%%%%%%%%%%%%%%%%%%%
IV) 
$b =0, \beta > 0$ & $ n^{-1} \sigma^{-(2a+1)}$ & $\sigma  >  n^{-\frac{1}{2k+2a+1}}$ & $ 0$ \\
                  & $n^{-\frac{2k}{2k+2a+1}}$ & $\sigma \leq n ^{-\frac{1}{2k +2a +1}}$ & $ n ^{-\frac{1}{2k +2a +1}}$\\ 
\hline

%%%%%%%%%%%%%%%%%%%%%
V) 
$\beta > b > 0$ & $ n^{-1} \exp\left(\kappa \sigma^{ \frac{-\beta b}{\beta - b}}\right)\sigma^ {\frac{\beta (b-2)}{2(\beta -b)}  -2\alpha }$ & 
$\sigma > \mu_1 $ & $ 0$ \\
 & $(\ln n)^{-\frac{2k}{b}}$ & $ \sigma \leq \mu_1 $ & $ \mu_1 $ \\
%%%%%%%%%%%%%%%%%%%%%%%
\hline
VI)
$b = \beta > 0$ & $ \sigma^{ -2\alpha} \left(\ln n \right)^{-\frac{2\alpha +2k}{b}}\, 
\exp\left(-2\gamma\sigma^\beta \left(\ln n \right)^{\frac{\beta}{b}}\right)$ & $\sigma > \mu_1 $ & $\mu_1  $ \\
 & $\left(\ln n \right)^{-\frac{2k}{b}}   $ & $ \sigma \leq \mu_1$ & $ \mu_2$ \\
%$\left(\frac{\ln n}{2(d- \gamma \sigma^{b}) }  \right)^{-\frac{2k}{b}}   $
 %%%%%%%%%%%%%%%%%%%%%%%%%%
\hline 
  VII)
  $b > 0,\beta = 0$ & $  \left(\ln n \right)^{-\frac{(2\alpha +2k)}{b}} \sigma^{-2\alpha} $ & $ \sigma  > \mu_1$ & $ \mu_1$ \\
 & $ \left(\ln n \right)^{-\frac{2k}{b}}  $ & $ \sigma  \leq  \mu_1$ & $ \mu_1$ \\
 %%%%%%%%%%%%%%%%%%%%% 
 \hline
 VIII)
  $b > \beta >  0$ & $  \sigma^{-2\alpha} \left(\ln n \right)^{\frac{(1+2a -2\alpha )}{b} -1}  $ & $ \sigma  >   \mu_1$ & $ \mu_1$ \\
 & $ \left(\ln n \right)^{-\frac{2k}{b}}  $ & $ \sigma  \leq  \mu_1$ & $ \mu_1$ \\
 \hline
 
\end{tabular}
%%%%%%%%%%%%%%%%%%%%%%%%%%%%%%%%%%%%%%%%%%%%%%%%%%%%%%%%%%%%%%%%%%%%%%%%%%%%%%%%%%%%%%%%
\caption{The   optimal values $h_{opt}$ of the bandwidth $h$ and the corresponding expressions for  $\Del (n,\sig, h)$ defined in \eqref{eq:Delta}.
 Here, $\mu_1$ and $\mu_2$ are given by \eqref{eq:mu}. 
}
\end{center}
 \end{table}

%%%%%%%%%%%%%%%%%%%%%%%%%%%%%%%%%%%%%%%%%%%%%%%%%%%%%%%%%%%%%%%%%%%%%%%%%%%%%%%%%%%%%%%%%

\section{Discussion}
\label{sec:discuss}
\setcounter{equation}{0}

In the present paper, our main goal was to theoretically justify the choice of a bandwidth in deconvolution problems 
with small Berkson errors. In particular, we  refined  the conclusion of Long \etal (2016) and studied the relationship between 
the three parameters: the bandwidth $h$, the sample size $n$ and the standard deviation of the Berkson errors $\sig$.
As Theorem~\ref{th:mise} above shows, the expressions for the optimal bandwidth is always chosen to minimize the error in the estimator 
of the density of interest $f_W$. In particular, if $h=0$ is possible, one should choose this value as long as the Berkson errors 
are not too small, i.e., $\sig$ lies above some threshold level that depends on the shapes of the densities and the number of 
observations $n$.

In order to uncover the reason for this, compare expressions \eqref{estimator} and  \eqref{ker_est} and observe that $g^* (\sigma s)$ in 
\eqref{estimator}  acts as a kernel function $g$ with  the bandwidth $h = \sigma$. If $\sigma$ is large enough 
(i.e $\sigma > h_{opt}$, where $h_{opt}$ is the value of $h$ that achieves the best bias-variance balance), 
then convolution with $g$ leads to sufficient regularization and no kernel estimation is necessary.
However , if $\sigma < h _{opt}$ then one needs additional kernel regularization with $h > \sigma $.
  
The setting of Long \etal (2016) corresponds to the cases I, II, III and IV in Tables 1 and 2 with $a=b=0$.
In all those cases, $h$ can be zero if $\sig$ is large enough or of the order $n^{-1/(2k+1)}$ where $k$ is the degree of smoothness of 
the   density $f_X$ of the measurements. The choice depends on the relationship between parameters $\sig$, $n$ and $k$.

Note that we did not consider the case of the multivariate density functions.
This extension is fairly straightforward but rather cumbersome. 
We shall leave this case for the future investigation.
 
%%%%%%%%%%%%%%%%%%%%%%%%%%%%%%%%%%%%%%%%%%%%%%%%%%%%%%%%%%%%%%%%%%%%%%%%%%%%%%%%%%%%%%%%%

\section*{Acknowledgments}

 Marianna Pensky  and Ramchandra Rimal were  partially supported by National Science Foundation
 (NSF), grants     DMS-1407475 and DMS-1712977.

%%%%%%%%%%%%%%%%%%%%%%%%%%%%%%%%%%%%%%%%%%%%%%%%%%%%%%%%%%%%%%%%%%%%%%%%%%%%%%%%%%%%%%%%%

% This is the version with reduced proofs.
% The complete proofs for the dissertation are in the file of Sept. 13, 2018

% This is the file with proofs

\section{Proofs}
\label{sec:proofs}
\setcounter{equation}{0}

\subsection{Proofs of the statements in the paper}
\label{subsec:main}

\noindent{\bf Proof of Lemma \ref{lbias}.\ }  Since 
\begin{align*}
 \Delta_1 & =\frac{1}{2\pi } \int_{| s| >1/ h}| g^*\left(\sigma s \right)| ^2 | f_X ^*  (s)| ^2 ds 
  =\frac{1}{\pi } \int_{\frac{1}{h}}^ {\infty} | g^*(\sigma s )| ^2 | f_X ^*  (s)| ^2 ds\\
& \leq \frac{2C_g}{\pi}\, \int_{\frac{1}{h}} ^{\infty} (\sigma^2 s^2 +1)^{-\alpha} \exp (-2\gamma| s|^{\beta} \sigma^{\beta} )
\frac{(s^2+1)^k }{(s^2 +1)^k}| f_X ^*  (s)| ^2 ds\\
&\leq \frac{2C_g\, B^2}{\pi}\,   \underset{s\geq\frac{1}{h}} {\max} \left[(\sigma^2 s^2 +1)^{-\alpha} \exp (-2\gamma| s|^{\beta} \sigma^{\beta} ) \right]\, 
(h^{-2} +1)^{-k},
\end{align*}
obtain
$$
\Delta_1 \asymp \min\left\{ \left(\frac{h}{\sigma}\right)^{2\alpha},1\right\} h^{2k} \exp\left(-2\gamma \left(\frac{\sigma}{h}\right)^{\beta}\right)
$$
which implies \eqref{bias}.
\\

\vspace{3mm}
 
%%%%%%%%%%%%%%%%%%%%%%%%%%%%%%%%%%%%%%%%%%%%%%%%%%%%%%%%%%%%%%%%%%%%%%%%%%%%%%%%%%%%%%%%%%%%%%%%%%%%%%%%%%%%%%%%%%%%%%
%%%%%%%%%%%%%%%%%%%%%%%%%%%%%%%%%%%%%%%%%%%%%%%%%%%%%%%%%%%%%%%%%%%%%%%%%%%%%%%%%%%%%%%%%%%%%%%%%%%%%%%%%%%%%%%%%%%%%%

\noindent
{\bf Proof of Lemma \ref{variance}.\ } 
Note that the  variance term is given by
\begin{align*}
\Delta_2   
& \leq \frac{1}{2\,\pi} \int_{-\infty}^{\infty}   \frac{  | g^* (\sigma s) | ^2} { | f_\xi ^*  (s)| ^2} I(| s| < h^{-1})  ds \\
&\leq \frac{C_g}{   c_\xi}\int_{0} ^\frac{1}{h} (\sigma^2 s^2 +1)^{-\alpha}\, (s^2 +1)^a \, 
\exp (-2\gamma| s|^{\beta} \sigma^{\beta}  +2d |s|^b)\, ds
\end{align*}
Using change of variables s=z/h obtain
\be\label{eq:Del2}     
 \Delta_2 \lesssim    h^{-(2a+1)}\, I(\sigma, h) \quad \mbox{with} \quad
I(\sigma, h) = \int_0^1 P(z| \sigma , h) \exp \{\phi(z | \sigma , h)\} dz
\ee
 where 
\be \label{phi_P}
 \phi(z | \sigma , h) =  2dz^b h^{-b} - 2\gamma z^{\beta}\sigma^{\beta} h^{-\beta}, \quad
P(z|\sigma,h) =    \left(\sigma^2  z^2 h^{-2} +1\right)^{-\alpha}\left( z^2  + h^{2}\right)^a
\ee
  
Let $d > 0, b > 0$. Denote by $z_0$ and $z_h$, respectively,  the global maximum of $\phi(z | \sigma , h)$ on the interval $[0,1]$ 
and its critical point
\be \label{eq:z0zh}
z_0 \equiv z_0(\sig, h) = \underset{z\in[0,1]} {\text{argmax}}\  \phi(z|\sigma,h), \quad 
z_h = \left( db\, (\gamma \beta)^{-1}\, \sigma^{-\beta} \right)^{\frac{1}{\beta-b}} h
% z_h = \left(\frac {\gamma \beta}{db}{\sigma^{\beta}}\right)^{\frac{1}{b - \beta}} h
\ee 
Then, since $z_h > 0$, $z_0 = z_h$ if $z_h \in (0,1]$ and $\phi''(z_h) < 0$ and $z_0=1$  otherwise.
Hence, Lemma \ref{lem:saddle_pt} and \eqref{eq:Del2}   yield  that  for small $h$ and    $\sigma $ 
\be \label{asymp_variance}
 h^{2a+1}\, \Delta_2 \lesssim \left\{
\begin{array}{ll}
   \frac{\exp\{ \phi(z_h|\sigma,h) \}P(z_h|\sigma,h)} {{\sqrt{ | \phi''(z_h|\sigma,h)|)}}}, &   \text{if} \quad  z_0 =z_h, \\
  \frac{\exp\{ \phi(1|\sigma,h) \}P(1|\sigma,h)} {{  \phi'(1|\sigma,h)}}, &   \text{if} \quad  z_0 = 1. 
\end{array}\right.
\ee 
Here
\be \label{expressions}
\begin{array}{ll}
 \phi(1 | \sigma , h) =  2d h^{-b} - 2\gamma \sigma^{\beta} h^{-\beta}, &  
{\phi}' (1| \sigma , h) = 2  \left( db h^{-b}  -  \gamma \beta \sigma^{\beta} h^{-\beta} \right) \\
P(1|\sigma,h) \asymp \left(\sigma^2   h^{-2} +1\right)^{-\alpha}, &  
P(z_h|\sigma,h) =    \left(\sigma^2  z_h^2 h^{-2} +1\right)^{-\alpha}\left( z_h^2  + h^{2}\right)^a
\end{array}
\ee
Below we consider  various cases.\\

%%%%%%%%%%%%%%%%%%%%%%%%%%%%%%%%%%%%%%%%%%%%%%%%%%%%%%%%%%

\noindent   \textbf{Cases I, II, III:} \     $b = \beta  =0. $ \\
Note that 
 \begin{align} \label{Del2CasesI-III}
\Delta_2 &\lesssim  \int_0^\frac{1}{h} (\sigma^2 s^2 +1)^{-\alpha}(s^2 +1)^a  ds % \\
 \asymp  h^{-1}\, \int_0^1 \left(\sigma^2 z^2 h^{-2} +1\right)^{-\alpha}\, \left( z^2 h^{-2} +1\right)^a\, dz 
 \end{align} 
If $h\geq  \sigma$, then $ \sigma^2 z^2 h^{-2}  + 1 \in (1,2)$ and  
$\Delta_2  \asymp   h^{-(2a+1)}$.
If  $h  < \sigma$,  then, by the change of variables $\sigma s = u$ in \eqref{Del2CasesI-III},   obtain 
 \begin{align*}
 \Delta_2 & \asymp  \sigma^{-1}\, \int_0^{\frac{\sigma}{h} } \left({u^2} +1\right)^{-\alpha}\,
\left( u^2 \sigma^{-2} +1 \right)^a\ du 
% \\& 
\asymp      {\sigma^{-(2a+1)}} \int_0^{\frac{\sigma}{h} }  \frac {u^{2a}} {{(u^2} +1)^{\alpha}}   du 
 \end{align*}
Hence,  
$$ 
\Delta_2 \asymp    \min \left( h^{-(2a +1)},  \sigma^{-(2a+1)}\right) \Delta_{h\sigma}
$$ 
where 
\be \label{Del_h_sigma}  
\Delta_{h\sigma} = \left\{
\begin {array}{lll}
1\    &  \text{if}\               \ \alpha > a + 1/2\\
  \max \left\{\ln\left({\frac{\sigma}{h} }\right),1\right\}\ & \text{if}\    \alpha = a + 1/2\\
\max\left\{ 1,\left({\frac{\sigma}{h} }\right)^{2a-2\alpha +1}\right \} & \text{if}\   \alpha < a + 1/2
\end {array}\right.
\ee 
\\

%%%%%%%%%%%%%%%%%%%%%%%%%%%%%%%%%%%%%%

\noindent   \textbf{Case IV:} \  $ b =0, \beta > 0 $. \\  
In this case, 
\bes %\label{var case 3}
 \Delta_2 \asymp  h^{-1}\, \int_0^1 \left(\sigma^2  z^2 h^{-2} +1 \right)^{-\alpha}\,
\left(z^2 h^{-2} +1\right)^a \exp \left(-2\gamma  \sigma^{\beta} z^{\beta}  h^{-\beta}   \right)\, dz.
 \ees
% Here $ \phi_h(z) = -2\gamma \sigma^{\beta} z^{\beta}  h^{-\beta}$.
If $h> \sigma$ then the argument of the exponent is bounded above and 
$\Delta_2 \asymp    h^{-2a-1}. $ 
If $h < \sigma$, then by changing variables  $u= 2 \gamma \left( \sigma z/h \right)^\beta,$ obtain
\bes
\Delta_2 \asymp  \sigma^{-1}\, \int_0^{\infty} 
\left(\left( \frac{u}{2\gamma } \right)^{\frac{2}{\beta}} +1\right)^{-\alpha}
\left(\frac{1}{\sigma^{2a}} \left(\frac{u}{2\gamma } \right)^{\frac{2a}{\beta}}+1\right)
 \exp (-u) u^{\frac{1}{\beta} -1} du \asymp    \sigma^{-(2a+1)}
\ees
  Hence,
 \bes 
 \Delta_2 \asymp   \min \left( h^{-(2a +1)},  \sigma^{-(2a+1)}\right).
\ees   
\\

%%%%%%%%%%%%%%%%%%%%%%%%%%%%%%%%%%%%%%

\noindent   \textbf{Case  V:} \  $\beta > b > 0$. \\  
In this case $\rho^2(\sigma) = \infty$ in \eqref{rho_delta}, so that  $h >0$.
The expression for the variance is given by \eqref{eq:Del2} with $ \phi(z | \sigma , h)$ defined in 
\eqref{phi_P}. Let  $z_h$ be given by \eqref{eq:z0zh}. It is easy to check that 
\be\label{opti_z}
  z_h = \left( db\, (\gamma \beta)^{-1}\, \sigma^{-\beta} \right)^{\frac{1}{\beta-b}} h \asymp \sigma^{-\frac{\beta}{\beta-b}}h.
\ee
It is easy to check that  ${\phi}''(z_h| \sigma , h) < 0$, so that $z_h$ is the local maximum. 
Now consider two cases.

(a) If $ h> \left(\frac {\gamma \beta}{db}{\sigma^{\beta}}\right)^{\frac{1}{\beta-b}}$,  then   
$z_h > 1$. Hence, $\phi(z|\sigma,h)$ does not have a local maximum on $[0,1]$ and it attains its global maximum at $z_0 = 1$.
Then,  
$2d h^{-b} >  \phi(1| \sigma , h) = 2d h^{-b} - 2\gamma \sigma^{\beta} h^{-\beta}  > 
2d h^{-b} \left(1 - b/\beta \right)$. 
Moreover, since $\beta > b$  and   $h> \left(\frac {\gamma \beta}{db}{\sigma^{\beta}}\right)^{\frac{1}{\beta-b}} > \sigma $, 
one has  $2dbh^{-b} > 2\gamma\beta \sigma^\beta h ^{-\beta}$ which yields 
$$
\phi^{'}(1| \sigma , h)=2db h^{-b} - 2\gamma\beta \sigma^\beta h ^{-\beta} =  
2dbh^{-b}\left(1 - \frac {\gamma \beta}{db}{\sigma^{\beta}}h^{b -\beta}\right)  \asymp h^{-b}
$$
Plugging those expressions into the second equation of  \eqref{asymp_variance} and using \eqref{expressions}, obtain
$$
\Delta_2 \asymp   h ^{ -(2a +1)} \min\left\{ \left( h \sigma^{-1}\right)^{2\alpha},1 \right\}\, 
\exp( 2dbh^{-b})\, h^b  
\asymp   h ^{b -2a -1} \exp\left(2dh^{-b}\right)
$$

(b) If $ h < \left(\frac {\gamma \beta}{db}{\sigma^{\beta}}\right)^{\frac{1}{\beta-b}}$,
then $z_h$ is given by formula \eqref{opti_z} and $z_0 = z_h <1$. 
Hence, $\Del_2$ is given by the first expression in formula \eqref{asymp_variance} 
\be  \label{Del2_caseVb} 
 \Delta_2 \asymp \frac{\exp(\phi(z_h| \sigma , h))}{\sqrt {|{\phi}''(z_h| \sigma , h)|}}\,  
h ^{ -(2a +1)}\left(\sigma^2  z_h^2 h^{-2} +1\right)^{-\alpha}\, \left(z_h^2 + h^{2}  \right)^a 
\ee 
Note that, due to  $\beta > b >0$, $\frac{\beta^2}{\beta-b} > \frac{\beta b}{\beta - b}$ 
and $\beta-\frac{\beta^2}{\beta-b} = -\frac{\beta b}{\beta - b}$, one has
\bes 
\phi(z_h| \sigma , h) = \frac{2d}{h^b}\left(\frac {db}{\gamma \beta}{\sigma^{-\beta}}\right)^{\frac{b}{\beta-b}}{h^b}  
-\frac{2\gamma \sigma^{\beta}}{h^{\beta}} \left(\frac {db}{\gamma \beta}{\sigma^{-\beta}}\right)^{\frac{\beta}{\beta-b}}{h^\beta}
= \kappa\sigma^{- \frac{\beta b}{\beta - b}}
\ees
where $\kappa$ is a positive constant defined in \eqref{kappa}. 
Also 
 $$
{\phi}'' (z_h| \sigma , h) = \frac{2}{z_h^2} \left(\frac{db(b-1) z_h^{b}}{h^{b}} -\frac{\gamma \beta(\beta-1)\ z_h^{\beta}\sigma^{\beta}} {h^{\beta}}\right)
 = \frac{2db(b-\beta)z_h^{b-2} }{h^b} \asymp \frac{z_h^{b-2} }{h^b}
$$ 
Then,   plugging $\phi(z_h| \sigma , h)$ and ${\phi}'' (z_h| \sigma , h)$ into \eqref{Del2_caseVb}, obtain 
\bes 
\Delta_2 \asymp \exp\left(\kappa  \sigma^{ -\frac{\beta b}{\beta - b}}\right)\sigma^{\frac{\beta (b-2)}{2(\beta -b)}  -2\alpha}.  
\ees
\\

%%%%%%%%%%%%%%%%%%%%%%%%%%%%%%%%%%%%%%

\noindent   \textbf{Case  VI:} \  $ b = \beta > 0$. \\ 
In this case $\rho^2(\sigma) = \infty$ in \eqref{rho_delta}, so that  $h >0$.
Moreover, since $\phi(z| \sigma , h)  =  2z^b\, h^{-b} (d- \gamma \sigma^{b})$ where, due to condition \eqref{sig-ineq},
$d- \gamma \sigma^{b}>0$, $z_0=1$ is the non-local maximum of $\phi(z| \sigma , h)$.  
 Then, the second expression in formula \eqref{asymp_variance}    
 \be \label{Del2_case2} 
\Delta_2 \lesssim \frac{\exp(\phi(1| \sigma , h))}{\phi' (1| \sigma , h)}\,  
h^{-(2a +1)}\left(\sigma^2  h^{-2} +1\right)^{-\alpha}    
\ee  
Using \eqref{expressions} with $\beta=b$, we derive
\bes
\Delta_2 \lesssim   h^{b-(2a+1)} \min\left(  \left( \frac{h}{\sigma} \right)^{2\alpha},1\right) \exp(2h^{-b}(d- \gamma \sigma^{b}))
\ees
\\

%%%%%%%%%%%%%%%%%%%%%%%%%%%%%%%%%%%%%%

\noindent   \textbf{Case  VII:} \  $ b >0, \beta = \gamma = 0 $\\
In this  case, $z_0=1$ is the non-local maximum of $\phi(z| \sigma , h)$  
and \eqref{expressions} yield $\phi(1 | \sigma , h) =  2d h^{-b}$ and 
$\phi' (1| \sigma , h) = 2   db h^{-b}$. Plugging those expressions into 
\eqref{Del2_case2}, we derive
\bes 
\Delta_2 \lesssim   \min\left(  \left(   \frac{h}{\sigma}  \right)^{2\alpha}  , 1  \right) h^{b-(2a+1)} {\exp (2dh^{-b}}) 
\ees
\\

%%%%%%%%%%%%%%%%%%%%%%%%%%%%%%%%%%%%%%

\noindent   \textbf{Case  VIII:} \ $ b > \beta > 0$\\
In this case $\rho^2(\sigma) = \infty$ in \eqref{rho_delta}, so that  $h >0$.
Also, it is easy to check that although $z_h \in (0,1)$, one has 
$\phi^{''}(z_h| \sigma , h) > 0$ , so $z_h$ is the local minimum.
It is easy to see that $z_0 =1$ and 
$\phi(1| \sigma , h) = {2d}{h^ {-b} }(1 - \gamma d^{-1}  \sigma^{\beta} h^{b- \beta}) \asymp  2d{h^ {-b} }.$
Moreover, 
$\phi' (1| \sigma , h) = 2     h^{-b} (db  -  \gamma \beta \sigma^{\beta} h^{b-\beta}) \asymp h^{-b}$,
so formula \eqref{Del2_case2} yields
\bes
\Delta_2 \lesssim   h^{b-(2a+1)} \min\left(  \left( \frac{h}{\sigma} \right)^{2\alpha},1\right) \exp(2dh^{-b}).
\ees
\\

\medskip

%%%%%%%%%%%%%%%%%%%%%%%%%%%%%%%%%%%%%%%%%%%%%%%%%%%%%%%%%%%%%%%%%%%%%%%%%%%%%%%%%
%%%%%%%%%%%%%%%%%%%%%%%%%%%%%%%%%%%%%%%%%%%%%%%%%%%%%%%%%%%%%%%%%%%%%%%%%%%%%%%%%%%%%%%%%%%%
\noindent
{\bf  Proof of Theorem \ref{th:mise}.\ } 
Consider various cases. 
\\

\noindent   
\textbf{Cases I, II, III:} \ \   $b = \beta  =0$.  \\
One has
\be \label{eq:Del_cases123}
\Delta \lesssim   \min\left\{ \left( h\, \sigma^{-1} \right)^{2\alpha}, 1\right\} h^{2k} + 
n^{-1} \min \left( h^{-(2a +1)},  \sigma^{-(2a+1)}\right) \Delta_{h\sigma}
\ee
where $\Delta_{h\sigma}$ is defined in \eqref{Del_h_sigma}. 
\\

%%%%%%%%%%%%%%%%%%%%%%%%%%%%%%%%%%%%%%%%%%%
%%%%%%%%%%%%%%%%%%%%%%%%%%%%%%%%%%%%%%%%%%%

\noindent  
\textbf {Case I:\ } $b = \beta  =0$,  $\alpha > a + 1/2$.  \\
In this case $\rho^2(\sigma) < \infty $ and $h =0$ is possible.  
If $h=0$, then   $\Delta = O \left({\sigma^{-(2a+1)}}{n^{-1}} \right).$ 
If $h \neq 0$, then  choose $h \geq \sigma$, so that $\Delta_1(\sigma , h) \lesssim h^{2k},$
$\Delta_2(\sigma , h) \lesssim h^{-(2a+1)}$. Then,  
$h_{opt}\asymp n^{-\frac{1}{2k+2a+1}}$ and 
$\Delta_1(\sigma , h_{opt}) + n^{-1}\, \Delta_2(\sigma , h_{opt}) \lesssim n^{-\frac{2k}{2k+2a+1}}$. 
Choose  $h = h_{opt}$ if $h_{opt} \geq \sigma$, i.e., if   $n^{-\frac{1}{2k+2a+1}}\geq \sigma$. 
Obtain 
\bes % \label{est case1a}
 \Delta \asymp  \left\{
\begin{array}{lll}
 n^{-1}\, \sigma^{-(2a+1)}, & h_{opt} =0  & \text{if}\  \sigma \geq n^{-\frac{1}{2k+2a+1}}\\
n^{-\frac{2k}{2k+2a+1}}, & h_{opt} =  n^{-\frac{1}{2k+2a+1}} & \text{if} \   \sigma  <  n^{-\frac{1}{2k+2a+1}}
\end{array} \right.
 \ees

%%%%%%%%%%%%%%%%%%%%%%%%%%%%%%%%%%%%%%%%%%%%%%%%%%

\noindent  \textbf{Case II:\ }   $b = \beta  =0$,  $\alpha = a +\frac{1}{2}.$\\
Here, $\Del$ is given by \eqref{eq:Del_cases123} where $\Delta_{h\sigma} = \max \lfi \ln(\sig/h), 1 \rfi$.
If $h<\sigma$, then 
$\Delta  \lesssim \sigma^{-2\alpha}h^{2\alpha+2k}+ \sigma^{-(2a+1)} n^{-1} \ln(\sig/h)$.
Setting   
$ \sigma^{-2\alpha} h^{2\alpha+2k} = \sigma^{-(2a+1)} n^{-1}\ln(\sig/h)$
leads to 
 $$
 h_{opt} \asymp   n^{-\frac{1}{2k+2\alpha}} \sigma^{\frac{2\alpha -2a-1}{ 2\alpha +2k}}, \quad 
 \Delta \lesssim  n^{-1}\, \sigma^{-(2a+1)}\,  \ln n.
  $$
 Note that 
 $ h_{opt} < \sigma$   if and only if $n^{-\frac{1}{2k+2a+1}}< \sigma$. 
Now, consider the  case when  $h\geq\sigma$. 
\noindent Then by \eqref{eq:Del_cases123}, 
$\Delta \lesssim  n^{-\frac{2k}{2k+2a+1}}$  if $n^{-\frac{1}{2k+2a+1}} \geq \sigma$. 
Hence 
\bes  % \label{est case1b}
 \Delta \asymp  \left\{
\begin{array}{lll}
\frac{\sigma^{-(2a+1)}}{n}\ln n, &  h_{opt} = n^{-\frac{1}{2k+2a}} \sigma^{\frac{2\alpha -2a-1}{ 2\alpha +2k}} &  
\text{if} \   \sigma  \geq  n^{-\frac{1}{2k+2a+1}}  \\
 n^{-\frac{2k}{2k+2a+1}}, & h_{opt}\asymp  n^{-\frac{1}{2k+2\alpha+1}}&  \text{if}\  \sigma < n^{-\frac{1}{2k+2a+1}}\  
\end{array} \right.
 \ees

%%%%%%%%%%%%%%%%%%%%%%%%%%%%%%%%%%%%%%%%%%%%%%%%%%

\noindent 
\textbf {Case III:\ }   $b = \beta  =0$,  $\alpha < a +\frac{1}{2}.$\\
First, consider the case when $h < \sigma$. Then, by   
\eqref{eq:Del_cases123} and   \eqref{Del_h_sigma}, obtain 
$$
\Delta \lesssim \sigma^{-2\alpha}h^{2\alpha+2k} +\sigma^{-(2\alpha)} n^{-1}  h^{2\alpha-2a -1}.  
$$
Setting
$ \sigma^{-2\alpha} h^{2\alpha+2k} = \sigma^{-(2\alpha)} n^{-1}  h ^{2\alpha-2a -1}$,
 obtain   $h_{opt} \asymp n^{-\frac{1}{2k +2a+1}}$ and 
$\Delta \lesssim \sigma ^{ -2\alpha} n ^{-\frac{2\alpha +2k}{2k +2a +1}}.$
Also note that $h_{opt} < \sigma$  if and only if $\sigma > n^{-\frac{1}{2k +2a+1}}$. 
Now, consider the case when $h\geq \sigma$. 
Then, \eqref{eq:Del_cases123} and   \eqref{Del_h_sigma}, derive that 
$\Delta  \asymp n^{-\frac{2k}{2k+2a+1}}$    if  $n^{-\frac{1}{2k+2a+1}} \geq \sigma$.
Hence 
\bes  % \label{est case1c}
 \Delta \asymp  \left\{
\begin{array}{lll}
 \sigma ^{-2\alpha } n ^{-\frac{2\alpha +2k}{2k +2a +1}}, & h_{opt} = n^{-\frac{1}{2k+2a+1}}   & \text{if}\   \sigma  >  n^{-\frac{1}{2k+2a+1}}  \\
n^{-\frac{2k}{2k+2a+1}}, & h_{opt} = n^{-\frac{1}{2k+2a+1}}  & \text{if} \    \sigma  \leq  n^{-\frac{1}{2k+2a+1}}    
\end{array} \right.
 \ees

%%%%%%%%%%%%%%%%%%%%%%%%%%%%%%%%%%%%%%%%%%%%%%%%%%

\noindent   
\textbf{Case  IV:\ }   $b =0, \beta > 0$.\\
In this case $\rho^2(\sigma) < \infty $ and $h =0$ is possible.  
Consider the case $h < \sigma$. Then,  
$$ 
\Delta_1(\sigma , h)   \lesssim  \sigma^{(-2\alpha)} h^{2\alpha +2k} \exp\left(-2\gamma \left(\frac{\sigma}{h} \right)^{\beta}\right),\quad
\Delta_2(\sigma , h)  \lesssim   \sigma^{-(2a+1)}.
$$ 
If $h < \sigma$,   then $h_{opt} =0$ and 
$\Delta \asymp  n^{-1} \sigma^{-(2a+1)}$. 
If  $h> \sigma$,  then
$\Delta_1(\sigma , h)  \leq   h^{2k}$ and  $\Delta_2(\sigma , h) \lesssim   h^{-(2a+1)}$. 
Therefore,   $h_{opt} \asymp n ^{-\frac{1}{2k +2a +1}}$ and 
$\Delta \lesssim n ^\frac{-2k}{2k +2a +1}$.
Observing that 
$h_{opt} \geq \sigma$ if $\sigma \leq n ^{-\frac{1}{2k +2a +1}}$, obtain
\bes
\Delta\asymp  \left\{
\begin{array}{lll}
 n^{-1} \sigma^{-(2a+1)} & h_{opt}  = 0 & \text{if}\   \sigma > n ^{-\frac{1}{2k +2a +1}} \\
  n ^{-\frac{2k}{2k +2a +1}} &	h_{opt} = n^{-\frac{1}{2k +2a +1}} & \text{if}\   \sigma \leq n ^{-\frac{1}{2k +2a +1}}
\end{array}\right.
\ees

%%%%%%%%%%%%%%%%%%%%%%%%%%%%%%%%%%%%%%%%%%%%%%%%%%

 \noindent   
\textbf{Case  V:\ }     $\beta > b > 0$.\\
In this case $\rho^2(\sigma) < \infty $ and $h =0$ is possible.  
The bias is given by \eqref{bias} and 
\bes
\Delta_2 \lesssim \left\{
\begin {array}{ll}
 n ^{-1} h ^{b -2a -1} \exp\left(2dh^{-b}\right) & \text{if}\ h > \left(\frac {\gamma \beta}{db}{\sigma^{\beta}}\right)^{\frac{1}{\beta-b}}\\
n^{-1} \exp\left(\kappa   \sigma^{ \frac{-\beta b}{\beta - b}}\right)\sigma^{\frac{\beta}{\beta -b}.\frac{b-2}{2} -2\alpha}  
& \text{if} \   h < \left(\frac {\gamma \beta}{db}{\sigma^{\beta}}\right)^{\frac{1}{\beta-b}}
\end {array}\right.
\ees
If $ h = 0 $, then
$ \Delta \asymp n^{-1} \exp\left(\kappa \sigma^{ \frac{-\beta b}{\beta - b}}\right)\sigma^{\frac{\beta}{\beta -b}.\frac{b-2}{2} -2\alpha}$.
If $h > 0$, then one needs $h > \sig \gtrsim \left(\frac {\gamma \beta}{db}{\sigma^{\beta}}\right)^{\frac{1}{\beta-b}}$ and 
$ \Delta \asymp  h^{2k}  + n ^{-1} h ^{b -2a -1} \exp\left(2dh^{-b}\right)$.
Choosing $h$ such that   $ h^{2k}  = n ^{-1} h ^{b -2a -1} \exp\left(2dh^{-b}\right)$,
arrive at 
 \be\label{exp_eqn_A}
(2dh^{-b})^\frac{2a+2k +1-b}{b} \exp\left(2dh^{-b}\right) =  (2d)^\frac{2a+2k +1-b}{b}\, n 
 \ee
and, by Lemma~\ref{lem:expoApprox}, 
obtain $ h_{opt}  =\mu_1(n)$ where $\mu_1(n)$ is defined in \eqref{eq:mu}, 
 and, hence,  $\Delta \asymp (\ln n)^{-\frac{2k}{b}}. $  
 Therefore, 
\bes
\Delta \asymp \left\{
 \begin {array}{lll}
 n^{-1} \exp\left(\kappa \sigma^{ \frac{-\beta b}{\beta - b}}\right)\sigma^  {\frac{\beta (b-2)}{2(\beta -b)}  -2\alpha }, 
& h_{opt} =0 & \text{if}\	\sigma \geq  \mu_1(n) \\
(\ln n)^{-\frac{2k}{b}}, & h _{opt}= \mu_1(n), &  
\text{if}\	 \sigma  < \mu_1(n)
\end{array}\right.
\ees
where $\mu_1(n)$ is given by \eqref{eq:mu}. 
\\

%%%%%%%%%%%%%%%%%%%%%%%%%%%%%%%%%%%%%%%%%%%%%%%%%%

\noindent   \textbf{Case  VI:} \  $b = \beta > 0, h>0$\\
Note that, due to \eqref{sig-ineq},    one has $\sigma < ( d \gamma^{-1})^\frac{1}{b}$.  
Consider two cases.  
If $h < \sigma$, then 
$$
\Delta_1(\sigma , h)  \lesssim   \sigma^{-2\alpha} h^{2\alpha +2k} \exp\left(-2\gamma \left( \sigma/h \right)^{\beta}\right), \quad  
\Delta_2(\sigma , h) \lesssim   h^{(b+2\alpha -2a-1)}  \sigma^{-2\alpha} \exp(2h^{-b}(d- \gamma \sigma^{b})). 
$$
Then the bias-variance balance is achieved when 
$$  
h^{(b-2k -2a-1)}   \exp(2h^{-b}(d- \gamma \sigma^{b}) + 2\gamma\sigma^b h^{ -b}  ) =  n
$$
which leads to \eqref{exp_eqn_A} and, hence,  $h_{opt} = \mu_1(n)$ where  $\mu_1(n)$ is defined in  \eqref{eq:mu}. 
Therefore, $h_{opt} \asymp \left(\ln n \right)^{-\frac{1}{b}}$    
and hence 
$$
\Delta \lesssim \sigma^{-2\alpha} \left(\ln n \right)^{-\frac{2\alpha +2k}{b}}\exp\left(-2\gamma\sigma^\beta 
\left(\ln n \right)^{\frac{\beta}{b}}\right).
$$
If $h \geq \sigma$, then 
$\Delta   \lesssim  h^{2k} + n^{-1} h^{b-(2a+1)} \exp(2h^{-b}(d- \gamma \sigma^{b}))$ 
and the  bias-variance balance is achieved when 
$ h^{2k} \asymp n^{-1} h^{b-(2a+1)} \exp(2h^{-b}(d- \gamma \sigma^{b}))$. 
Then, by  Lemma \ref{lem:expoApprox}, we derive that 
$ h_{opt} = \mu_2(n)$ where  $\mu_2(n)$ is defined in  \eqref{eq:mu}, and 
$ \Delta  \lesssim (\ln n) ^{-\frac{2k}{b}}.$
Hence 
\bes
\Delta \lesssim \left\{
\begin{array}{lll}
 \sigma^{(-2\alpha)} \left(\ln n \right)^{-\frac{2\alpha +2k}{b}}\, \exp\left(-2\gamma\sigma^\beta \left(\ln n \right)^{\frac{\beta}{b}}\right),&
h _{opt}= \mu_1(n),  & \text{if}\  \sigma > \mu_1(n) \\
\left(\ln n  \right)^{-\frac{2k}{b}}, &
h _{opt}= \mu_2(n),  & \text{if}\  \sigma \leq \mu_1(n)  
\end{array}\right.
\ees
where $\mu_1(n)$ and $\mu_2(n)$ are given by \eqref{eq:mu}. 
\\

%%%%%%%%%%%%%%%%%%%%%%%%%%%%%%%%%%%%%%%%%%%%%%%%%%

\noindent   \textbf{Case  VII}:  \  $ b >0,  \beta = 0 $\\
If $h < \sigma$, then 
$$
\Delta \lesssim   \sigma ^{-2\alpha}h^{2\alpha + 2k} +   n^{-1}\, \sigma^{-2\alpha}\,  h ^{2\alpha  -2a +b -1} \exp(2dh^{-b})
$$
Setting 
$\sigma ^{-2\alpha}h^{2\alpha + 2k} =   n^{-1}\, \sigma^{-2\alpha}\,  h ^{2\alpha  -2a +b -1} \exp(2dh^{-b})$,
arrive at \eqref{exp_eqn_A} and $h_{opt} = \mu_1(n)$ where  $\mu_1(n)$is defined in  \eqref{eq:mu}. 
Hence, 
$h_{opt} \asymp \left(\ln n \right)^{-1/b}$ and 
$ 
\Delta \lesssim    \left(\ln n \right)^{-\frac{ 2\alpha +2k}{b}}\,  \sigma^{-2\alpha}, 
$ 
provided $\sig > \mu_1(n)$.

If $h \geq \sigma$, then 
\be\label{mse case 2}
 \Delta \lesssim h^{ 2k} + n^{-1} h ^{b -2a -1} {\exp(2dh^{-b})}. 
\ee
Setting $h^{2k} \approx  n^{-1} h ^{b -2a -1}  \exp(2dh^{-b})$,  
arrive at \eqref{exp_eqn_A}, so that $h_{opt} = \mu_1(n) \asymp \left(\ln n \right)^{-1/b}$
and $\Delta \lesssim   (\ln n)^{ -(2k/b}$ if $\sig \leq \mu_1(n)$.
Hence 
\bes  % \label{est case2}
 \Delta \asymp  \left\{
\begin{array}{lll}
 \left(\ln n \right)^{-\frac{2\alpha +2k}{b}} \sigma^{-2\alpha}, &  h _{opt}=  \mu_1(n), &	\text{if}\	 \sigma >   \mu_1(n) \\
 \left(\ln n \right)^{-\frac{2k}{b}}, &  h _{opt}=  \mu_1(n),  & \text{if}   \sigma  \leq  \mu_1(n),
\end{array} \right.
 \ees
where  $\mu_1(n)$is defined in  \eqref{eq:mu}. 
\\

%%%%%%%%%%%%%%%%%%%%%%%%%%%%%%%%%%%%%%%%%%%%%%%%%%%%%%%%%%%%%%%%%%%%%%%%%%%%%%%%%%%%%%%%

 \noindent   \textbf{Case  VIII}: \ $ b > \beta > 0$  \\
If $h \leq \sigma$, then 
$$
\Delta(\sigma,h)\lesssim \sigma^{-2\alpha} h^{2\alpha +2k} \exp\left(-2\gamma {\sigma}^\beta {h}^{-\beta}\right)  + 
n^{-1} h^{2\alpha + b-(2a+1)}  {\sigma} ^{-2\alpha} \exp(2dh^{-b}).
$$
Then, the  minimum of $\Delta(\sigma , h)$ is attained if 
$n \asymp h^{b-(2a+1)-2k}   \exp(2dh^{-b} +2\gamma {\sigma}^\beta {h}^{-\beta})$.
Note that, due to ${\sigma}^\beta < (d/\gamma)\, h ^{-(b-\beta)}$,  $b > \beta$ and $\sigma < 1$, one has
 $2dh^{-b} > 2\gamma {\sigma}^\beta {h}^{-\beta}$ . Therefore, we arrive at \eqref{exp_eqn_A},  so that 
 $h_{opt} \asymp \left(\ln n \right)^{- 1/b}$ and  
$\Delta \lesssim \sigma^{ -2\alpha}  (\ln n)^{\frac{(1+2a -2\alpha )}{b} -1}$.

If  $h> \sigma$, then 
$\Delta \lesssim   h^{2k} + n^{-1} h^{b-(2a+1)}  \exp(2dh^{-b})$
which coincides with \eqref{mse case 2} and we obtain the same expressions for $h_{opt}$ and $\Delta$ as in that case.  
Hence 
\bes  % \label{est case2}
 \Delta \asymp  \left\{
\begin{array}{lll}
 \sigma^{-2\alpha} \left(\ln n \right)^{\frac{(1+2a -2\alpha )}{b} -1} , &  h _{opt}=  \mu_1(n), &	\text{if}\	 \sigma >   \mu_1(n) \\
 \left(\ln n \right)^{-\frac{2k}{b}}, &  h _{opt}=  \mu_1(n),  & \text{if}   \sigma  \leq  \mu_1(n),
\end{array} \right.
 \ees
where  $\mu_1(n)$is defined in  \eqref{eq:mu}. 
\\
 
 % 
%\bes
%\Delta\asymp  \left\{
%\begin{array}{ll}
%  \sigma^{-2\alpha} \left(\ln n \right)^{\frac{(1+2a -2\alpha )}{b} -1} ,h _{opt}= \left[\frac{1}{2d}\left( \ln(A) + \frac{b -2a-2k -1}{b}\ln\ln(A)\right)\right]^{-\frac{1}{b}}\\ 
%  \quad 	 \text{if}\	   \sigma > \left[ \frac{1}{2d}\left(\ln n + \frac{b-2k-2a-1}{b}\ln\ln n\right)\right]^{-\frac{1}{b}}  \\
%  
% \left(\ln n \right)^{-\frac{2k}{b}} ,h _{opt}= \left[\frac{1}{2d}\left( \ln(A) + \frac{b -2a-2k -1}{b}\ln\ln(A)\right)\right]^{-\frac{1}{b}}\\
% \quad   \text{if}\  \sigma \leq \left[ \frac{1}{2d}\left(\ln n + \frac{b-2k-2a-1}{b}\ln\ln n\right)\right]^{-\frac{1}{b}}    
%\end{array}\right.
%\ees 

%%%%%%%%%%%%%%%%%%%%%%%%%%%%%%%%%%%%%%%%%%%%%%%%%%%%%%%%%%%%%%%%%%%%%%%%%%%%%%%%%%%%%%%%%%%%%%%%%%%%%%%%%%%%%%%%%%%%%%%%%%%%%%%%%%%%%%%
%%%%%%%%%%%%%%%%%%%%%%%%%%%%%%%%%%%%%%%%%%%%%%%%%%%%%%%%%%%%%%%%%%%%%%%%%%%%%%%%%%%%%%%%%%%%%%%%%%%%%%%%%%%%%%%%%%%%%%%%%%%%%%%%%%%%%%%

\subsection{Supplementary statements and their proofs}
\label{subsec:suppl}

\begin{lemma}\label{lem:saddle_pt}
Consider an integral of the form
\be \label{int_lem1}
I(\lambda) = \displaystyle \int_{m_1}^{m_2} P_{\lambda}(z) \exp(Q_{\lambda}(z)) dz
\ee 
where $0 \leq m_1 < m_2 < \infty$ and $P_{\lambda}(z)$ and $Q_{\lambda}(z)$ are real valued differentiable functions of $z$
and  $\lambda \to \infty$ is a large parameter.
Let 
\bes
 z_0 \equiv z_{0,\lambda}   = \underset{z\in[m_1,m_2]} {\text{argmax}}Q_{\lambda}(z)
\ees 
be an unique  global maximum of $Q_{\lambda}(z)$ on the interval $[m_1, m_2]$. Assume that the following conditions hold:
\begin{itemize}

\item 
$P_{\lambda}(z)$ is a slow varying function that can be expanded into Taylor series at $z = z_0$
and such that $|P_{\lambda}(z)| \leq M$  for some constant $M$ independent of $\lambda$;

\item
$Q_{\lambda}(z_0) - Q_{\lambda}(z)$  increases steadily  to $\infty$ as $\lambda \to \infty$  

\item
If $Q_{\lambda}^{'} (z_0) = 0$, then for every $ \lambda \geq \lambda_0$  
\be  \label{eq:local}
\lim_{x \to 0}  \frac{Q_{\lambda}(z_0 +x)  - Q_{\lambda}(z_0)}{x^2} = \frac{Q_{\lambda}^{''}(z_0)}{2} < 0
\ee

\item
If $Q_{\lambda}^{'} (z_0) \neq 0$, then for every $ \lambda \geq \lambda_0$ 
\be  \label{eq:nonlocal}
\lim_{x \to 0}  \frac{Q_{\lambda}(z_0 + x)  - Q_{\lambda}(z_0)}{x} =  Q_{\lambda}^{'}(z_0)  \neq 0
\ee

\end{itemize}

\noindent
Then, as $\lambda \to \infty$, 
\be \label{eq:saddle_pt}
I(\lambda) \ \asymp    \lfi 
\begin{array}{ll}
\frac{\exp\{ Q_{\lambda}(z_0) \}\, P_{\lambda}(z_0)}{  {\sqrt{ | Q_{\lambda}^{''}(z_0)|}}}, & \mbox{if  \quad \eqref{eq:local} holds},\\
\frac{\exp\{Q_{\lambda}(z_0) \}\, P_{\lambda}(z_0)}{ Q_{\lambda}^{'}(z_0)}, & \mbox{if  \quad \eqref{eq:nonlocal} holds}
\end{array} \right.
\ee
\end{lemma}

\noindent
Validity of Lemma \ref{lem:saddle_pt} follows from Sections 5.2 and 5.3 of Dingle (1973).
In particular, the first expression in  \eqref{eq:saddle_pt} follows from formulas (3) and (4), page 111,
 and the second one follows from formulas (16) and (17), page 119 of Dingle (1973).

\medskip
\vspace{4mm}

\begin{lemma} \label{lem:expoApprox}
Let $n$ be large. Then solution of the equation 
\be \label{exp_approx}
e^m m^z = n
\ee 
is given by  
\be \label{exp_solution}
m \approx \ln n - z \ln \ln n.
\ee 
\end{lemma}
\medskip

\noindent
{\textbf{Proof}. } 
Since  $e^m m^z = n$, then  $ m + z \ln m = \ln n$ and  $m = \ln n -z \ln m$.
Plugging this $m$ back into \eqref{exp_approx}, obtain $e^{\ln n -z \ln m } (\ln n -z \ln m )^z = n$.
Since  for large values of $n$, one has $(\ln n -z \ln m )^z \approx (\ln n)^z$, the previous equation becomes 
$(\ln n)^z  n e^{ -z \ln m }  \approx n$, so that $z \ln \ln n  \approx  z \ln m$ which yields \eqref{exp_solution}. \\

%%%%%%%%%%%%%%%%%%%%%%%%%%%%%%%%%%%%%%%%%%%%%%%%%%%%%%%%%%%%%%%%%%%%%%%%%%%%%%%%%%%%%%%%%%%%%%%%%%%%%%%%%%%%%%%%%%%%%%%%%%%%%%

%%%%%%%%%%%%%%%%%%%%%%%%%%%%%%%%%%%%%%%%%%%%%%%%%%%%%%%%%%%%%%%%%%%%%%%%%%%%%%%%%%%%%%%%%%%%
%%%%%%%%%%%%%%%%%%%%%%%%%%%%%%%%%%%%%%%%%%%%%%%%%%%%%%%%%%%%%%%%%%%%%%%%%%
 
% \newpage

\end{document}